\title{Unsuitable Mathematics: From Card Games to Probability}
\author{Tiernan Colby, Baojun Fang, Mushtaq Kapasi, Yanle Li, \and Cameron Miller, Kayla Miller, Steven J. Miller, Yilong Wu}
\date{}
\documentclass{article}
\usepackage{graphicx}
\usepackage{amsmath,amsfonts,amsthm,amssymb,amscd}
\usepackage{color}
\usepackage{hyperref}
\usepackage{url}
\newcommand{\bburl}[1]{\textcolor{blue}{\url{#1}}}
\usepackage{float}
\usepackage{tikz}
\usepackage{tkz-tab}
\usepackage{tkz-graph}
\usetikzlibrary{shapes.geometric,positioning}

\newcommand{\p}{\mathbb{P}}
\newtheorem{theorem}{Theorem}

\newtheorem{definition}{Definition}
\newtheorem{corollary}{Corollary}
\newtheorem{remark}{Remark}

\renewenvironment{proof}{{\sc Proof:}}{~\hfill QED}

\begin{document}
\newpage
\maketitle
\begin{abstract}
We introduce a two-player card dealing game in which each player receives one card per round. A player is considered suited if all cards in that player's hand belong to the same suit and unsuited otherwise. The game ends when both players are unsuited or when the entire deck of cards is exhausted, whichever occurs first. We determine the expected duration of the game. Along the way, we use the game as an accessible introduction to common tools in probability and combinatorics that are useful in research such as the inclusion-exclusion principle, tail-sum formulas for expectations, and limits of probability functions.  
\end{abstract}




\section{Introduction} 

Probability theory has its roots, in part, from analyzing games of chance, dating back at least to the seventeenth century. Since then, the study of card games has long provided a rich source of problems and inspired new directions in probability research. For example, Pierre Rémond de Montmort's investigation of a card-matching game in the early eighteenth century led to the study of derangements \cite{Mon08}, while Bayer and Diaconis's analysis of card shuffling in 1992 provided important insights into the convergence of random processes to equilibrium \cite{BD92}. These examples illustrate how seemingly simple questions can give rise to substantial mathematical investigations.
Here, we use a simple card-dealing game as an accessible introduction to tools that are used in probability and combinatorics research. These tools include the inclusion-exclusion principle, the tail-sum formula for expected values, conditional probability, recurrences, and limits to show convergence. By analyzing the game, we hope readers will obtain a deeper understanding of these tools being used and will learn techniques that they can apply to their own problems.

The card game we analyze in this paper originated from a dinner conversation that Steven Miller had with his son Cam and his daughter Kayla on low-probability risk and gambling strategies. Their discussion eventually turned to a simple card-dealing game and led to the following question: On average, how many rounds of cards must be dealt before every player holds cards from at least two distinct suits?

Motivated by this question, we formulate the following game. Two players are dealt cards from a deck of $N$ cards divided equally among $s$ suits, where $s \ge 2$ and $N \geq 4$ is an even integer.\footnote{The assumption that $N$ is even ensures that both players receive the same number of cards if the entire deck is dealt. This makes our analysis of the game easier, because both players would have an equal number of cards by the end of every round. We require $N\ge 4$ to exclude the trivial case $N=2$, in which the game necessarily ends after the first round.} Each player receives one card each round, and a round is considered complete only after both players have received their cards. We consider that a player is suited if all her cards belong to the same suit and unsuited if her cards belong to more than one suit. The game ends when either both players are \emph{unsuited} or all cards have been dealt, whichever occurs first. Let the random variable $T$ denote the round which the game ends. Our primary goal is to compute $\mathbb{E}[T]$, which is the expected value of $T$. $\mathbb{E}[T]$ represents the average round at which the game ends if the game was played many times.

\begin{figure}[htbp]
    \centering
    \includegraphics[width=0.9\textwidth]{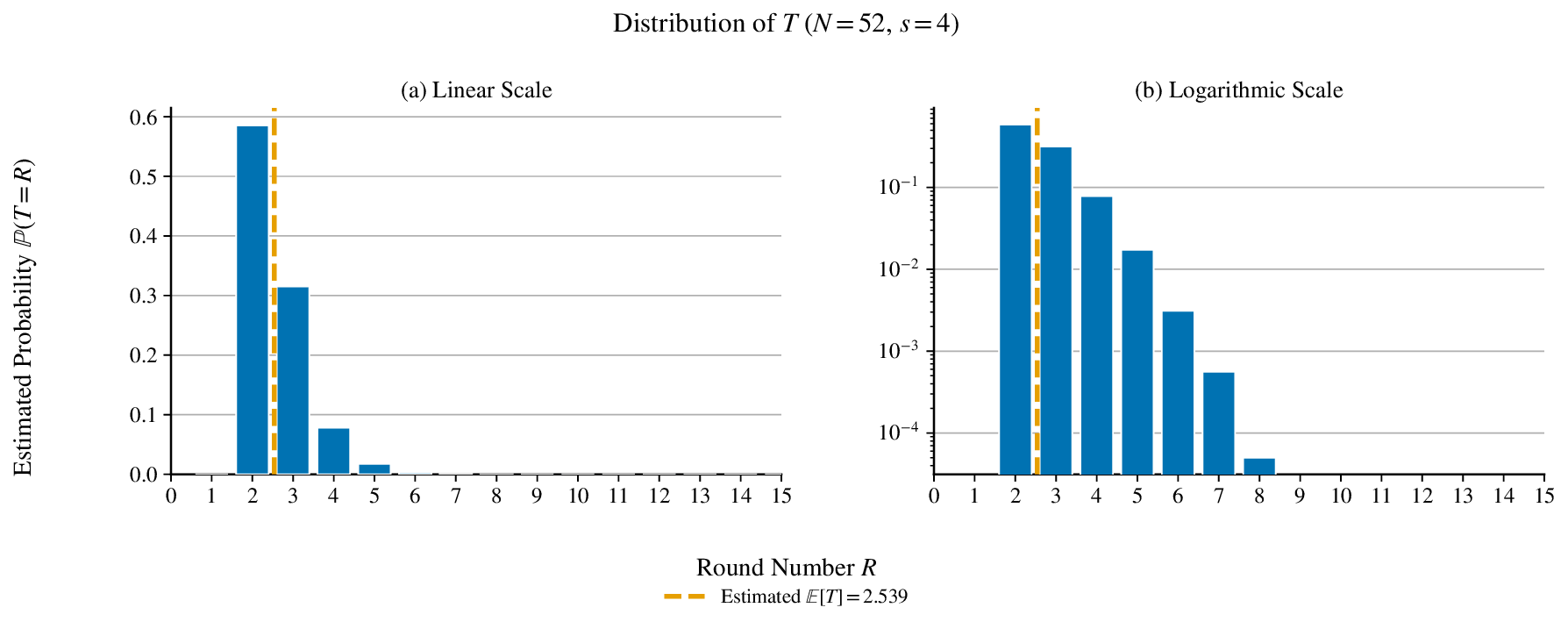}
    \caption{Estimated distribution of $T$, the round that the game ends, from Monte Carlo simulations using a standard deck of cards $(N=52, s=4)$ and $100{,}000$ trials.} \label{fig:simulation_card_dealing}
\end{figure}

As a first step to tackle the problem, we conducted simulations to estimate the expected number of rounds until the game ends. We simulated the game using a standard deck of cards $(N = 52, s = 4)$, and the result is shown in Figure \ref{fig:simulation_card_dealing}. 
Based on the simulation result, we see that the estimated $\mathbb{E}[T]$ is around $2.539$, which is the yellow dotted line in the figure. We observe that most of the probability mass is concentrated in early rounds, and the probability decays rapidly as the rounds get larger. As the probabilities in later rounds are extremely small, we also present the log graph to show the tail behavior.

Our plan for computing $\mathbb{E}[T]$ is to first derive $\mathbb{P}(T = R)$, the probability mass function (PMF) that the game ends at $R$, and then use it to compute the expected value. Because the game ends when both players become unsuited, we aim to calculate the PMF by examining the probability that both players are unsuited by a given round $R$. The main obstacle to our plan, however, is the dependence between successive draws. In particular, the cards dealt in the previous rounds affect the probabilities of drawing cards of a given suit in the subsequent rounds, and the cards received by one player also influence the distribution of cards available to the other player. Thus, it is hard to derive the PMF directly. When encountering problems of this kind, we are reminded of George Pólya’s famous principle, which John Conway attributes to Pólya in the foreword to his book\emph{How to Solve It} \cite{Pol04}: ``\emph{if you can't solve a problem, then there is an easier problem you can solve: find it}''. Here, we see that the finite number of cards creates the obstacle, so we ask whether we can solve the simpler case of an infinite deck containing $s$ equally likely suits.

This paper is organized as follows. In Section \ref{sec:toolkit}, we define and explain the tools from probability and combinatorics that we use throughout the paper. We include definitions of the PMF, expected value, the inclusion-exclusion principle, and explain our strategy of computing the probability that the game ends at round $R$ as the difference of the probabilities that it has ended by rounds $R$ and $R-1$.

In Section \ref{sec:infinite-two-suit}, we study this infinite deck version of the game. The key simplification is that the probability of drawing a card of a given suit no longer depends on the cards dealt in the previous rounds, and the distribution of cards available to each player at each round is now the same. Thus, it is possible to compute $\mathbb{P}(T = R)$, the probability mass function that the game ends at a fixed round $R$, by multiplying the probabilities that Player 1 and Player 2 are unsuited at round $R$ and then differencing consecutive rounds as in \eqref{eq:ptr-pmf-from-cdf}. We can then compute the expected value via the PMF and obtain  
\[
\mathbb E[T] \ = \ \frac{2s^2+2s-1}{s^2-1}.
\]
Because the probability of drawing a card of a given suit does not change from round to round, we can also compute $\mathbb{E}[T]$ with a recurrence in Section~\ref{subsec:infinite-round-recurrence} without finding the PMF. 

In Section \ref{sec:finite-deck}, we return to the case of the finite deck. Unlike the infinite case, we no longer have independence. To overcome this difficulty, we first consider the easier case that the game has not ended by round $R$. Define $Q_R = \p(T > R)$ for $1 \le R \le N /2$ as the probability that the game has not ended by round $R$. To compute $\mathbb{E}{T}$, we use the inclusion-exclusion principle to determine the probability that the game has not ended. We take a combinatorial approach by counting the number of hands that result in the game continuing. Using this, we solve for the expected value of $T$: 
\[
\mathbb{E}[T] \ = \ 2+\sum_{R=2}^{N/2-1} Q_R, 
\]
where 
\[ 
Q_R
\ = \ 2s\,\frac{\binom{N/s}{R}}{\binom{N}{R}}
- s\,\frac{\binom{N/s}{2R}}{\binom{N}{2R}}
- s(s-1)\,\frac{\binom{N/s}{R}^2}{\binom{N}{R}\binom{N-R}{R}}. 
\]

\section{Problem-Solving Toolkit} 
\label{sec:toolkit}
While the game is simple to formulate, its analysis requires a lot of classical tools and concepts encountered in probability and combinatorics. Here, we present valuable tools used in this paper. For readers interested in learning about those tools in more detail, \cite{Mil17} is a great reference. 

\begin{definition}[Probability Mass Function]
The probability mass function (PMF) of a discrete random variable $T$ is
\[ 
f_T(t) \ = \ \mathbb{P}(T = t),
\]
the probability that $T$ takes exactly the value $t$.
\end{definition}

\begin{definition}[Expected Value] 
The expected value of a discrete random variable $T$ is 
\[ 
\mathbb{E}[T] \ = \ \sum_{t} t\, \mathbb{P}(T = t).  
\]
It can be interpreted as the weighted average of the possible values of $T$, where each value is weighted by its probability.
\end{definition}
In our game, the random variable $T$ represents the round at which the game ends. Thus, $\mathbb{E}[T]$ describes the average number of rounds required for the game to end over many independent repetitions. One natural approach to computing this quantity is to first determine the PMF of $T$ and then apply the definition of expected value.
\begin{definition}[Inclusion-Exclusion Principle]
For any two events $A_1$ and $A_2$, we have
\[ 
\mathbb{P}(A_1 \cup A_2) \ = \ \mathbb{P}(A_1) + \mathbb{P}(A_2) - \mathbb{P}(A_1 \cap A_2).
\]
\end{definition}
This principle will be particularly useful in our finite-deck analysis, where the game continues if at least one player remains suited. Another useful strategy is to compute the probability that a random variable takes a value of at most $t$, rather than \emph{exactly} the value $t$. 
\begin{definition}[Cumulative Distribution Function] 
The cumulative distribution function (CDF) of a discrete random variable $T$, denoted $F_T(t)$, is 
\[ 
F_{T}(t) \ = \ \mathbb{P}(T \leq t), 
\]
which gives the probability that $T$ takes a value at most $t$.
\end{definition}

Throughout the paper, we compute the probability that the game ends at round $R$ as the probability that the game has ended by round $R$ minus the probability that the game has ended by round $R-1$. In terms of the CDF, this probability is
\begin{equation}
\mathbb{P}(T = R) \ = \ \mathbb{P}(T \leq R) - \mathbb{P}(T\leq R- 1). 
\label{eq:ptr-pmf-from-cdf}
\end{equation}
Equivalently, expressed as the probability of the game having not ended yet,
\begin{equation}
\mathbb{P}(T = R) \ = \ \mathbb{P}(T > R-1) - \mathbb{P}(T > R).
\label{eq:ptr-tail-difference}
\end{equation}
Throughout the rest of the paper, we use $\mathbb{P}(\text{ends by } R)$ to denote the CDF $\mathbb{P}(T \leq R)$.  We frequently use this strategy because it is often easier to calculate the probability of the game ending by round $R$ than the probability of the game ending exactly at round $R$. 

\section{Warm-up: Infinite Deck}
\label{sec:infinite-two-suit}

\subsection{The Probability That One Player Is Suited}
\label{subsec:infinite-one-player}
For a fixed suit $j\in\{1,\dots,s\}$, we say that a player is $j$-suited if all the player's cards have suit $j$. At round $R\geq1$, Player $1$ is $j$-suited exactly when all $R$ cards have suit $j$. Because every suit is equally likely, we have 
\begin{equation}
\p(\text{Player 1 is }j\text{-suited at }R)
\ = \ \left(\frac1s\right)^R.
\label{eq:infinite-single-player-fixed-suit}
\end{equation}
We note that $\{\text{Player 1 is suited at } R\} = \bigsqcup_{j = 1}^{s} \{\text{Player 1 is } j\text{-suited at } R\}$, where the union is disjoint because a player cannot be suited in more than one suit at the same time. This then gives 
\begin{align}
\p(\text{Player 1 suited at }R)
&\ = \ \sum_{j=1}^{s}\p(\text{Player 1 is }j\text{-suited at }R)
\nonumber\\
&\ = \ s\left(\frac1s\right)^R \nonumber \\
&\ = \ s^{1-R}.
\label{eq:infinite-single-player-matching}
\end{align}
Taking the complement, we obtain
\begin{equation}
\p(\text{Player 1 unsuited at }R) \ = \ 1-s^{1-R}, 
\label{eq:infinite-single-player-unsuited}
\end{equation}
for $R \geq 1$. As a natural check, we see that $\p(\text{Player 1 unsuited at }R) = 0$ when $R = 1$, and it increases as $R$ increases.  

\subsection{The PMF and the Expected Value}
\label{subsec:infinite-two-player}

Once a player becomes unsuited, adding more cards cannot make the player suited again. Thus, the game has ended by round $R$ exactly when both players are unsuited at that round. The two players' card sequences are independent, so we multiply the individual probabilities in \eqref{eq:infinite-single-player-unsuited} to obtain the following distribution.

\begin{theorem} 
\label{thm:infinite-two-suit}
In the infinite deck game, the probability mass function is given by 
\begin{equation}
\p(T=R) \ = \ 
\begin{cases}
\displaystyle
\frac{2(s-1)}{s^{R-1}}-\frac{s^2-1}{s^{2R-2}},
& R\geq2,\\[6pt]
0,&R = 1. 
\end{cases}
\label{eq:infinite-game-pmf}
\end{equation}
\end{theorem}

\begin{proof}
As remarked, at $R = 1$, each player has only one card and is suited, so $\mathbb{P}(T = 1) = 0$. We then compute the probability that the game has ended by round
$R$ for $R \geq 1$. Independence and
\eqref{eq:infinite-single-player-unsuited} give
\begin{align}
\p(\text{ends by }R)
&\ = \ \p(\text{Player 1 unsuited at } R \, \cap \, \text{Player 2 unsuited at } R)
\nonumber\\
&\ = \ \p(\text{Player 1 unsuited at }R) 
  \, \p(\text{Player 2 unsuited at }R)
\nonumber\\
&\ = \ \left(1-s^{1-R}\right)^2.
\label{eq:infinite-game-cdf-proof}
\end{align}
To find the probability that the game ends exactly at round $R$, we use \eqref{eq:ptr-pmf-from-cdf}. We subtract the probability that it has ended by round $R-1$ from the probability that it has ended by round $R$:
\begin{align}
\p(T=R)
&\ = \ \p(\text{ends by }R)-\p(\text{ends by }R-1)
\nonumber\\
&\ = \ \left(1-s^{1-R}\right)^2-\left(1-s^{2-R}\right)^2
\nonumber\\
&\ = \ 2(s-1)s^{1-R}-(s^2-1)s^{2-2R},
\label{eq:infinite-game-pmf-proof}
\end{align}
which proves \eqref{eq:infinite-game-pmf} for $R\geq2$.  
\end{proof}

\begin{corollary} 
\label{cor:infinite-game-moments}
Under the assumptions of Theorem~\ref{thm:infinite-two-suit},
\begin{equation}
\mathbb{E}[T] \ = \ \frac{2s^2+2s-1}{s^2-1}.
\label{eq:infinite-first-moment-result}
\end{equation}
\end{corollary}

\begin{proof}
We use the standard power-series identity
\begin{equation}
\sum_{i=1}^{\infty}i x^i \ = \ \frac{x}{(1-x)^2},
\qquad |x|<1, 
\label{eq:infinite-series-first-moment}
\end{equation}
which can be easily obtained by differentiating the geometric series and multiplying
by $x$; see \cite[Eq.~(1.2.2)]{Wilf94} or \cite{Mil17}.

Let $i\geq1$ count the additional rounds after the initial deal,
so the corresponding stopping round is $i+1$. Because the game ends
with probability one, $\sum_{i=1}^{\infty}\p(T=i+1)=1$.
Separating the initial round from the additional rounds and
substituting \eqref{eq:infinite-game-pmf}, we obtain
\begin{align}
\mathbb{E}[T]
&\ = \ \sum_{i=1}^{\infty}(i+1)\p(T=i+1)
\nonumber\\
&\ = \ \sum_{i=1}^{\infty}\p(T=i+1)
  +\sum_{i=1}^{\infty}i\,\p(T=i+1)
\nonumber\\
&\ = \ 1+2(s-1)\sum_{i=1}^{\infty}\frac{i}{s^i}
     -(s^2-1)\sum_{i=1}^{\infty}\frac{i}{s^{2i}}
\nonumber\\
&\ = \ 1+\frac{2s}{s-1}-\frac{s^2}{s^2-1}
\nonumber\\
&\ = \ \frac{2s^2+2s-1}{s^2-1}.
\label{eq:infinite-first-moment}
\end{align}
The two series converge because $s\geq2$; their values follow from
\eqref{eq:infinite-series-first-moment} with $x=1/s$ and $x=1/s^2$.
\end{proof}

\subsection{A Recurrence for the Expected Value}
\label{subsec:infinite-round-recurrence}
In Section~\ref{subsec:infinite-two-player}, we computed $\mathbb{E}[T]$ by first finding the PMF of $T$. In the case of infinite deck, however, a recurrence approach allows us to compute $\mathbb{E}[T]$ without finding the PMF. The key observation is that we can view the game as moving from one state to another, where the state at the end of a round is the number of suited players: $0$, $1$, or $2$. 

For $k\in\{0,1,2\}$, let $E_k$ denote the expected number of additional rounds until the game ends, given that exactly $k$ players are suited at the end of the current round. If no player is suited, the game has already ended, so $E_0=0$ by definition. At $R=1$, each player has only one card and is suited, so the game is in state $2$ after the first round. Thus, $T$ consists of the first round together with the additional rounds played from state $2$, which gives
\begin{equation}
\mathbb{E}[T] \ = \ 1+E_2.
\label{eq:infinite-expected-state-two}
\end{equation}
Therefore our goal is to compute $E_2$.  
We first compute the probabilities of transitioning from one state to another and list them in Table \ref{tab:infinite-transition}. 

\begin{table}[H]
\centering
\renewcommand{\arraystretch}{1.5}
\begin{tabular}{c|ccc}
\hline
 & \multicolumn{3}{c}{Next state} \\
Current state & $0$ & $1$ & $2$ \\
\hline
$0$ & $1$ & $0$ & $0$ \\
$1$ & $\dfrac{s-1}{s}$ & $\dfrac{1}{s}$ & $0$ \\
$2$ & $\dfrac{(s-1)^2}{s^2}$ & $\dfrac{2(s-1)}{s^2}$ & $\dfrac{1}{s^2}$ \\[6.5pt]
\hline
\end{tabular}
\caption{Transition probabilities between states, where state $k$ means $k$ players are suited. Each row gives the current state and each column gives the state after the next round.}
\label{tab:infinite-transition}
\end{table}
\noindent Conditioning on the three outcomes of the next round, we obtain
\begin{align}
E_2
&\ = \ 1+\frac{1}{s^2}E_2+\frac{2(s-1)}{s^2}E_1+\frac{(s-1)^2}{s^2}E_0 \nonumber\\
&\ = \ 1+\frac{1}{s^2}E_2+\frac{2(s-1)}{s^2}E_1.
\label{eq:infinite-recurrence-E2}
\end{align}
We also obtain 
\begin{equation}
E_1 \ = \ 1+\frac{1}{s}E_1+\left(1-\frac{1}{s}\right)E_0 \ = \ 1+\frac{1}{s}E_1.
\label{eq:infinite-recurrence-E1}
\end{equation}
Solving for $E_1$ gives $E_1=s/(s-1)$.
Substituting $E_1=s/(s-1)$ into \eqref{eq:infinite-recurrence-E2}, the $E_1$ term becomes $2/s$ and collecting the terms involving $E_2$, we have
\begin{equation}
E_2 \ = \ \frac{s^2+2s}{s^2-1}.
\end{equation}
By \eqref{eq:infinite-expected-state-two}, we conclude that
\begin{equation*}
\mathbb{E}[T] \ = \ 1+\frac{s^2+2s}{s^2-1} \ = \ \frac{2s^2+2s-1}{s^2-1},
\end{equation*}
which is the same as \eqref{eq:infinite-first-moment-result}.

\section{Results on Finite Deck}\label{sec:finite-deck}
\subsection{Combinatorial Argument Overview} 
\label{sec:argument-overview}
We now return to the finite deck game. Unlike the game setup in Section \ref{sec:infinite-two-suit}, the probability of drawing a card is now dependent on the remaining cards because both players draw from a finite deck of size $N$. The independence argument and the recurrence argument in Section \ref{sec:infinite-two-suit} do not apply. Thus, we use a combinatorial strategy to derive the PMF $\mathbb{P}(T = R)$ and compute $\mathbb{E}[T]$.

To compute $\mathbb{P}(T = R)$, we first consider the two boundary cases $R = 1$ and $R = N/2$. At $R = 1$, both players have only one card and are suited, so the game cannot end. Thus, we have $\mathbb{P}(T = 1) = 0$. At $R = N/2$, all cards have been dealt, so the game must end. Note that the game ends at $R = N/2$ only if it has not ended at $R = N/2 - 1$, so $\mathbb{P}(T = N/2) = \mathbb{P}(T > N/2 - 1)$, which is \eqref{eq:ptr-tail-difference} with $\mathbb{P}(T > N/2) = 0$. To obtain $\mathbb{P}(T = R)$ for intermediate rounds $1 < R < N/2$, we first compute the probabilities that the game has not ended by round \(R-1\) and by round \(R\), and then subtract the latter from the former as in \eqref{eq:ptr-tail-difference}. Let $Q_R = \p(T > R)$ denote the probability that the game has not ended by round $R$. Because a player who becomes unsuited remains unsuited for the rest of the game, the game has not ended by round $R < N/2$ exactly when at least one player is suited at round $R$. By the inclusion-exclusion principle, we write $Q_R$ as  
\begin{align}
Q_R
&\ = \ \mathbb{P}\left( \text{Player 1 suited at } R \cup \text{Player 2 suited at } R \right) \nonumber\\
&\ = \ \mathbb{P}(\text{Player 1 suited at } R) + \mathbb{P}(\text{Player 2 suited at } R) \nonumber \\ 
&\qquad \qquad - \mathbb{P}(\text{Players 1 and 2 both suited at } R).
\label{eq:finite-inclusion-exclusion}
\end{align}

The deck is shuffled uniformly at random, so any set of $R$ cards is equally likely to be any other size $R$ set. Thus, each player's hand has the same distribution, so $\mathbb{P}(\text{Player 1 suited at } R)$ and $\mathbb{P}(\text{Player 2 suited at } R)$ are equal. Thus we can simplify the expression to
\begin{equation} 
Q_R \ = \  2\mathbb{P}(\text{Player 1 suited at } R) - \mathbb{P}(\text{Players 1 and 2 both suited at } R).
\label{eq:master}
\end{equation} 

The rest of this section is organized as follows. In Section~\ref{subsec:finite-probabilities}, we compute the probability that one player is suited at round $R$ and the probability that both players are suited at $R$. In Section~\ref{subsec:finite-pmf} we compute $Q_R$ and derive the PMF $\mathbb{P}(T = R)$. In Section~\ref{subsec:finite-expectation}, we use the PMF to calculate the average length of the game. In Section~\ref{subsec:finite-convergence} we show that the finite deck probabilities converge to the infinite deck probabilities as the deck size increases.


\subsection{Computation of Probabilities}\label{subsec:finite-probabilities} 
We first compute $\p(\text{Player 1 suited at $R$)}$. At the end of round $R$, the player has received $R$ cards from the deck of $N$ cards, so there are $\binom{N}{R}$ hands. There are $N/s$ cards of a fixed suit $j$ in the deck, so there are $\binom{N/s}{R}$ possible hands out of $\binom{N}{R}$ total hands. For any fixed suit $j$, we have 
\begin{equation}
\mathbb{P}(\text{Player } 1 \text{ is } j \text{-suited at } R)
\ = \ \frac{\binom{N/s}{R}}{\binom{N}{R}}, \qquad j = 1,\dots,s.
\label{eq:finite-single-player-fixed-suit}
\end{equation}
Summing over all $s$ suits, we obtain
\begin{equation}\label{eq:PSA}
\mathbb{P}(\text{Player 1 suited at } R) \ = \ s\,\frac{\binom{N/s}{R}}{\binom{N}{R}} .
\end{equation}

Next, we calculate $\mathbb{P}(\text{Players 1 and 2 both suited at } R)$. To calculate this probability, we consider the case where both players are suited in the same suit and when the players are suited in different suits. In the case where both players have the same suit at $R$, the $2R$ cards held by the players must come from the same suit. Similar to our argument for \eqref{eq:PSA}, we get 
\begin{equation}
\mathbb{P}(\text{Players 1 and 2 both suited in the same suit at } R)
\ = \ s\,\frac{\binom{N/s}{2R}}{\binom{N}{2R}}.
\label{eq:finite-both-same-suit}
\end{equation}

Now we consider the case where both players have different suits. We track which suit each player holds. Suppose that Player 1 is $j$-suited, and Player 2 is $k$-suited, for $j \neq k$. We see that Player 1 draws $R$ cards of suit $j$ with probability $\binom{N/s}{R}/\binom{N}{R}$. Given that Player 1 is $j$-suited, Player $2$ draws $R$ cards of suit $k$ from the $N - R$ cards with probability $\binom{N/s}{R}/\binom{N-R}{R}$. This then gives the following joint probability:  
\begin{align}
\mathbb{P}( \substack{\text{Player 1 } j\text{-suited and} \\   \text{Player 2 }   k \text{-suited at } R})
&\ = \ \mathbb{P}(\text{Player 1 } j\text{-suited at } R) \nonumber \\ &\qquad \times \mathbb{P}(\text{Player 2 }  k\text{-suited at } R\mid \text{Player 1 } j\text{-suited at } R) \nonumber  \\ 
&\ = \ \frac{\binom{N/s}{R}}{\binom{N}{R}} \frac{\binom{N/s}{R}}{\binom{N-R}{R}}.
\label{eq:finite-ordered-suit-assignment}
\end{align}
This is the probability for one ordered suit assignment, and because there are $s$ choices for player $1$'s suit and $s-1$ remaining choices for player $2$'s, there are $s(s-1)$ assignments, so
\begin{align}
\mathbb{P}(\text{Players 1 and 2 suited in different suits at } R)
&\ = \ s(s-1)\,\frac{\binom{N/s}{R}}{\binom{N}{R}}\frac{\binom{N/s}{R}}{\binom{N-R}{R}} \nonumber \\ 
&\ = \ s(s-1)\,\frac{\binom{N/s}{R}^2}{\binom{N}{R}\binom{N-R}{R}}.
\label{eq:finite-both-different-suits}
\end{align}
Adding these two cases together, we obtain 
\begin{equation}\label{eq:PSAB}
\mathbb{P}(\text{Players 1 and 2 both suited at } R)
\ = \ s\,\frac{\binom{N/s}{2R}}{\binom{N}{2R}}
+ s(s-1)\,\frac{\binom{N/s}{R}^2}{\binom{N}{R}\binom{N-R}{R}}.
\end{equation}

\begin{remark}
Both players can hold cards of the same suit only if $2R \le N/s$, because the $2R$ cards dealt by round $R$ must come from a suit containing $N/s$ cards. The formula in \eqref{eq:PSAB} is valid for larger $R$ because $\binom{a}{b} = 0$ when $b > a$, so the same-suit term goes to $0$ when that event is not possible. The terms involving $\binom{N/s}{R}$ are $0$ once $R > N/s$.
\end{remark}

\subsection{The PMF}
\label{subsec:finite-pmf} 
By substituting \eqref{eq:PSA} and \eqref{eq:PSAB} into \eqref{eq:master}, we get
\begin{equation}\label{eq:QR}
Q_R
\ = \ 2s\,\frac{\binom{N/s}{R}}{\binom{N}{R}}
- s\,\frac{\binom{N/s}{2R}}{\binom{N}{2R}}
- s(s-1)\,\frac{\binom{N/s}{R}^2}{\binom{N}{R}\binom{N-R}{R}} 
\end{equation}
for $1 \le R < N/2$.
By \eqref{eq:ptr-pmf-from-cdf}, the probability that the game ends exactly at round $R$ is  
\begin{align}
\p(T = R)
&\ = \ \p(\text{ends by } R) - \p(\text{ends by } R-1) \nonumber  \\
&\ = \ \bigl(1 - Q_R\bigr) - \bigl(1 - Q_{R-1}\bigr) \nonumber \\
&\ = \ Q_{R-1} - Q_R,   
\label{eq:QR-1-QR}
\end{align}
for $1 < R < N/2$.

Collecting \eqref{eq:QR-1-QR} and our observations about boundary cases in Section \ref{sec:argument-overview}, we have 
\begin{equation}\label{eq:finitepmf}
\mathbb{P}(T = R) \ = \ 
\begin{cases}
0 & R = 1, \\ 
Q_{R-1} - Q_R   & 2 \le R < N/2, \\
Q_{N/2-1} & R = N/2, 
\end{cases}
\end{equation}
where we see 
\begin{align}
Q_{R - 1} - Q_R
&\ = \ 2s\left(\frac{\binom{N/s}{R-1}}{\binom{N}{R-1}} - \frac{\binom{N/s}{R}}{\binom{N}{R}}\right)  - s\left(\frac{\binom{N/s}{2R-2}}{\binom{N}{2R-2}} - \frac{\binom{N/s}{2R}}{\binom{N}{2R}}\right) \nonumber \\
&\qquad  - s(s-1)\left(\frac{\binom{N/s}{R-1}^2}{\binom{N}{R-1}\binom{N-R+1}{R-1}} - \frac{\binom{N/s}{R}^2}{\binom{N}{R}\binom{N-R}{R}}\right),
\label{eq:finite-pmf-intermediate}
\end{align}
by substituting \eqref{eq:QR} into \eqref{eq:QR-1-QR}. 

\subsection{Computing the Expected Value}
\label{subsec:finite-expectation}
\begin{corollary} \label{cor:finite-exp}
The expected round at which the finite deck game ends is
 \begin{equation}\label{eq:finite-first-moment-explicit}
\mathbb{E}[T]
\ = \ 2+\sum_{R=2}^{N/2-1}Q_R,
\end{equation}
with $Q_R$ defined in \eqref{eq:QR}.

\end{corollary}
\begin{proof} 
For the discrete random variable $T$, the expected value is
\begin{equation}
\mathbb{E}[T] \ = \ \sum_{R} R \, \p(T = R).
\label{eq:finite-expectation-definition}
\end{equation}
By \eqref{eq:finitepmf}, for $T \in \{2,\dots,N/2\}$, $\p(T=R)=Q_{R-1}-Q_R$ for $2\le R<N/2$ and $\p(T=N/2)=Q_{N/2-1}$. Because the last round receives the remaining probability, we can split the sum into the final round and the rounds before it,
\begin{equation}
\mathbb{E}[T] \ = \ \sum_{R=2}^{N/2-1} R (Q_{R-1} - Q_{R}) + \frac{N}{2}Q_{N/2-1}.
\label{eq:finite-expectation-split}
\end{equation}
Writing out these terms, 
\begin{align}
\mathbb{E}[T]
&\ = \ 2(Q_1 - Q_2)  + \cdots + (N/2-1)(Q_{N/2-2} - Q_{N/2-1}) + (N/2)\,Q_{N/2-1} \nonumber \\
&\ = \ 2Q_1 + (3-2)Q_2 + (4-3)Q_3 + \cdots + \bigl(N/2-(N/2-1)\bigr)Q_{N/2-1} \nonumber \\
&\ = \ 2Q_1 + Q_2 + Q_3 + \cdots + Q_{N/2-1}.
\label{eq:finite-expectation-telescoping}
\end{align}
Each $Q_R$ term for $2 \le R \le N/2-1$ appears in the sum once with a coefficient of $R+1$ and once with a coefficient of $-R$, which sums to a coefficient of $1$. Recall that $Q_1 = 1$, so $2Q_1 = 2$. This gives \eqref{eq:finite-first-moment-explicit}.
\end{proof} 


\begin{remark} 
    In fact, it is possible to obtain $\mathbb{E}[T]$ without computing $\mathbb{P}(T = R)$. The tail sum formula $\mathbb{E}[T] = \sum_{R\ge 0}\p(T>R)$, which applies to random variables taking nonnegative integer values, gives Corollary \ref{cor:finite-exp}. This can also be used to calculate the infinite deck result in Corollary~\ref{cor:infinite-game-moments}. In both game setups, $\mathbb{E}[T] = 2 + \sum_{R\ge 2} Q_R$ and the two games differ in $Q_R$, the probability that the game has not ended. 
\end{remark}


\subsection{Convergence Analysis}
\label{subsec:finite-convergence}
We now show pointwise convergence of the finite deck PMF to the infinite deck PMF by letting $N$ tend to infinity through even integers divisible by $s$. We keep the number of suits $s\geq2$ and the round $R\geq2$ fixed. For sufficiently large $N$, we have $R<N/2$, so the formula for intermediate rounds in \eqref{eq:finitepmf} applies.

Recall from \eqref{eq:finite-pmf-intermediate} that  
\begin{align}
\mathbb{P}(T = R) &\ = \ 2s\left(\dfrac{\binom{N/s}{R-1}}{\binom{N}{R-1}} - \dfrac{\binom{N/s}{R}}{\binom{N}{R}}\right)
- s\left(\dfrac{\binom{N/s}{2R-2}}{\binom{N}{2R-2}} - \dfrac{\binom{N/s}{2R}}{\binom{N}{2R}}\right) \nonumber \\ 
&\qquad - s(s-1)\left(\frac{\binom{N/s}{R-1}^2}{\binom{N}{R-1}\binom{N-R+1}{R-1}} - \frac{\binom{N/s}{R}^2}{\binom{N}{R}\binom{N-R}{R}}\right), 
\label{eq:finite-convergence-pmf-intermediate}
\end{align}
for $2 \leq R < N/2$. 
We first rewrite each fraction of binomial coefficients as a finite product. Take the following ratio as an example; we can expand it and cancel the common factor $R!$:
\begin{align}
\frac{\binom{N/s}{R}}{\binom{N}{R}}
&\ = \ \frac{(N/s)(N/s-1)\cdots(N/s-R+1)}{N(N-1)\cdots(N-R+1)} \nonumber \\
&\ = \ \prod_{i=0}^{R-1}\frac{N/s-i}{N-i} \nonumber \\
&\ = \ \prod_{i=0}^{R-1}\frac{1/s-i/N}{1-i/N}.
\label{eq:finite-convergence-product}
\end{align}

Taking the limit, we see that 
\begin{equation}
\lim_{N\to\infty}\frac{\binom{N/s}{R}}{\binom{N}{R}}
\ = \ \prod_{i=0}^{R-1}
  \lim_{N\to\infty}\frac{1/s-i/N}{1-i/N}
\ = \ \frac{1}{s^R}. 
\label{eq:finite-convergence-single-limit}
\end{equation}
Replacing $R$ by $R-1$, $2R-2$, and $2R$ in the same calculation gives
\begin{equation}
\label{eq:finite-convergence-other-ratios}
\lim_{N\to\infty}\frac{\binom{N/s}{R-1}}{\binom{N}{R-1}}
\ = \ \frac{1}{s^{R-1}}, 
\quad \lim_{N\to\infty}\frac{\binom{N/s}{2R-2}}{\binom{N}{2R-2}}
\ = \ \frac{1}{s^{2R-2}}, 
\quad \lim_{N\to\infty}\frac{\binom{N/s}{2R}}{\binom{N}{2R}}
\ = \ \frac{1}{s^{2R}}.
\end{equation}
For the term in which the two players hold different suits, we separate the ratio into two factors:
\begin{equation}
\begin{aligned}
\frac{\binom{N/s}{R}^{2}}
     {\binom{N}{R}\binom{N-R}{R}}
&\ = \ \frac{\binom{N/s}{R}}{\binom{N}{R}}
   \frac{\binom{N/s}{R}}{\binom{N-R}{R}}. 
\end{aligned}
\label{eq:finite-convergence-different-product}
\end{equation}
As each factor converges, we have 
\begin{equation}
\lim_{N\to\infty}
\frac{\binom{N/s}{R}^{2}}
     {\binom{N}{R}\binom{N-R}{R}}
\ = \ \frac{1}{s^{2R}}.
\label{eq:finite-convergence-different-limit}
\end{equation}
By the same reasoning as above, we also have 
\begin{equation}
\lim_{N\to\infty}
\frac{\binom{N/s}{R-1}^{2}}
     {\binom{N}{R-1}\binom{N-R+1}{R-1}}
\ = \ \frac{1}{s^{2R-2}}.
\label{eq:finite-convergence-different-previous}
\end{equation}
Substituting these limits into \eqref{eq:finite-pmf-intermediate}, and combining the terms with coefficients $s$ and $s(s-1)$, we obtain
\begin{align}
\lim_{N\to\infty}\p(T=R)
&\ = \ 2s\left(\frac{1}{s^{R-1}}-\frac{1}{s^R}\right)
   -s\left(\frac{1}{s^{2R-2}}-\frac{1}{s^{2R}}\right)\nonumber \\
&\qquad
   -s(s-1)\left(\frac{1}{s^{2R-2}}-\frac{1}{s^{2R}}\right)\nonumber \\
&\ = \ \frac{2(s-1)}{s^{R-1}}
   -\frac{s^2-1}{s^{2R-2}}.
   \label{eq:finite-convergence-pmf-limit}
\end{align}
The probability on the left refers to the finite deck game with $N$ cards. The expression on the right is the infinite deck probability in \eqref{eq:infinite-game-pmf}.

\begin{figure}[htbp]
    \centering
    \includegraphics[width=0.75\linewidth]{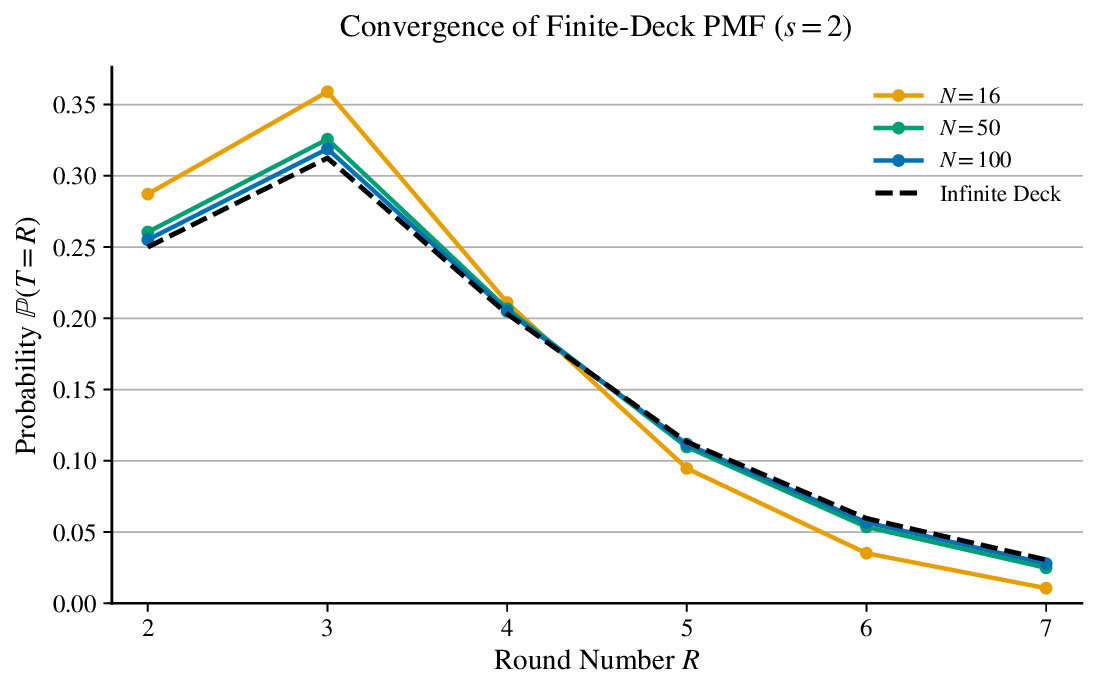}
    \caption{Infinite deck and Finite deck with varying deck size}
    \label{fig:finite-inf-convergence}
\end{figure}

To illustrate the convergence, we compare decks with two suits and $N=16$, $50$, and $100$ cards, with $N/2$ cards of each suit. We evaluate the finite deck probabilities \eqref{eq:finitepmf} for $2 \le R \le 7$ and compare them to the infinite deck probabilities \eqref{eq:infinite-game-pmf}. In Figure~\ref{fig:finite-inf-convergence}, the dashed line shows the infinite deck probabilities from \eqref{eq:infinite-game-pmf}. The blue line ($N=100$) is closest to the dashed line, the green line ($N=50$) is second closest, and then there is a gap to the orange line ($N=16$). The finite deck probabilities approach the infinite deck probabilities as $N$ increases.

\section{Future Work}\label{sec:conclusion}

We studied the average round at which a two-player card-dealing game ends. There are several natural extensions that follow from our work. First, we assumed that each suit has the same number of cards in the finite deck. However, if the suits contain different numbers of cards, the inclusion–exclusion formula still applies, and we could replace \(N/s\) by the number of cards in the chosen suit and sum over the probability of all the possible assignments of suits. Second, it is interesting to compute the results of the game in the case of $P$ players, and we believe our argument can still apply though with small complications. We leave these extensions to the interested reader.

\section{Acknowledgment}
This paper is partially supported by NSF Grant DMS2341670. We thank Mackenzie Pitsko for helpful discussions during the preparation of this paper. AI tools, including ChatGPT and Claude, were used to assist with the literature review, coding simulations, and proofreading of this manuscript.

 

\section*{About the authors:}

\subsection*{Tiernan Colby}
   Department of Mathematics and Statistics,
   Pomona College,
   Claremont, CA 91711, USA.
   teca2023@mymail.pomona.edu

\subsection*{Baojun Fang}
   Department of Mathematics,
   Stanford University,
   Stanford, CA 94305, USA.
   bjfang26@stanford.edu

\subsection*{Mushtaq Kapasi}
   Department of Mathematics,
   Open University,
   Milton Keynes, UK.
   mushtaq.kapasi@aya.yale.edu

\subsection*{Yanle Li}
   School of Mathematics and Physics,
   Xi'an Jiaotong-Liverpool University,
   Suzhou, 215123, PRC.
   yanle.li23@student.xjtlu.edu.cn

\subsection*{Cameron Miller}
   Florida Institute of Technology,
   Melbourne, FL 32901, USA.
   cmill17sport@gmail.com

\subsection*{Kayla Miller}
   Department of Mathematics,
   Williams College,
   Williamstown, MA 01267, USA.
   ksmill5433@gmail.com

\subsection*{Steven J. Miller}
   Department of Mathematics and Statistics,
   Williams College,
   Williamstown, MA 01267, USA.
   sjm1@williams.edu,
   Steven.Miller.MC.96@aya.yale.edu

\subsection*{Yilong Wu}
   Department of Mathematics and Statistics,
   Williams College,
   Williamstown, MA 01267, USA.
   mw17@williams.edu
\end{document}